\documentclass[aos,MSNbibl,nameyear,dvips]{arximspdf}

\doi{10.1214/12-AOS973} 
\volume{40}
\issue{1}
\pubyear{2012}
\firstpage{609}
\lastpage{637}

\makeatletter

\newtheorem{theorem}{Theorem}[section]
\newtheorem{proposition}[theorem]{Proposition}
\newtheorem{corollary}[theorem]{Corollary}
\newtheorem{lemma}[theorem]{Lemma}

\newproclaim{definition}[theorem]{Definition}
\newproclaim{example}[theorem]{Example}
\newproclaim{remark}[theorem]{Remark}

\newcommand{\cal}{\mathcal}

\newcommand{\real}{{\mathbb{R}}}
\newcommand{\Real}{\overline{\real}}

\newcommand{\cA}{{\cal A}}
\newcommand{\cD}{{\cal D}}

\newcommand{\cM}{{\cal M}}
\newcommand{\cP}{{\cal P}}
\newcommand{\cR}{{\cal R}}
\newcommand{\cS}{{\cal S}}

\newcommand{\myS}{{\mathrm S}}
\newcommand{\mys}{{\mathrm s}}

\newcommand{\pd}{\partial}

\makeatother

\begin{document}
\begin{frontmatter}

\title{Local proper scoring rules of order two}
\runtitle{Local proper scoring rules}

\begin{aug}
\author[A]{\fnms{Werner} \snm{Ehm}\ead[label=e1]{ehm@igpp.de}}
\and
\author[B]{\fnms{Tilmann} \snm{Gneiting}\corref{}\thanksref{t1}\ead[label=e2]{t.gneiting@uni-heidelberg.de}}
\runauthor{W. Ehm and T. Gneiting}
\affiliation{Institute for Frontier Areas of Psychology and Mental
Health and~University of Heidelberg}
\address[A]{Institute for Frontier Areas\\
\quad of Psychology \\
\quad and Mental Health \\
Wilhelmstr. 3a~\\
79098 Freiburg \\
Germany\\
\printead{e1}}
\address[B]{Institute for Applied Mathematics \\
University of Heidelberg \\
Im Neuenheimer Feld 294 \\
69120 Heidelberg \\
Germany\\
\printead{e2}} 
\end{aug}

\thankstext{t1}{Supported by the Alfried Krupp von Bohlen und Halbach
Foundation, and by the United States National Science Foundation
under Awards ATM-0724721 and DMS-07-06745.}

\received{\smonth{2} \syear{2011}}
\revised{\smonth{1} \syear{2012}}

%
\begin{abstract}
Scoring rules assess the quality of probabilistic forecasts, by
assigning a~numerical score based on the predictive distribution and
on the event or value that materializes. A~scoring rule is proper if
it encourages truthful reporting. It is local of order $k$ if the
score depends on the predictive density only through its value and the
values of its derivatives of order up to $k$ at the realizing event.
Complementing fundamental recent work by Parry, Dawid and Lauritzen,
we characterize the local proper scoring rules of order 2 relative
to a~broad class of Lebesgue densities on the real line, using a~different approach. In a~data example, we use local and nonlocal
proper scoring rules to assess statistically postprocessed ensemble
weather forecasts.
\end{abstract}

%

\pdfkeywords{62C99, 62M20, 86A10, Density forecast,
Euler equation, Hyvarinen score, proper scoring rule,
tangent construction}

\end{frontmatter}

\section{Introduction} \label{secintroduction}

One of the major purposes of statistical analysis is to make forecasts
for the future, and to provide suitable measures of the uncertainty
associated with them. Consequently, forecasts ought to be
probabilistic in nature, taking the form of probability distributions
over future quantities and events [\citet{Daw84}]. Scoring rules provide
summary measures for the evaluation of probabilistic forecasts, by
assigning a~numerical score based on the predictive distribution and
on the event or value that materializes. We take scoring rules to be
negatively oriented losses that a~forecaster wishes to minimize.
Specifically, if the forecaster quotes the predictive distribution $Q$
and the event $x$ materializes, her loss is $\myS(x,Q)$. The function~$\myS(\cdot,Q)$ takes values in the extended real line, $\Real=
[-\infty,\infty]$, and we write $\myS(P,Q)$ for the expected value of
$\myS(\cdot,Q)$ under $P$. Suppose, then, that the forecaster's
best judgment is the predictive distribution $P$. The forecaster has
no incentive to predict any $Q \not= P$, and is encouraged to quote
her true belief, $Q = P$, if $\myS(P,P) \leq\myS(P,Q)$. A~scoring
rule with this property is said to be proper [\citet{GneRaf07}].

Our paper is concerned with local proper scoring rules for
probabilistic forecasts of a~real-valued quantity. Briefly, if the
predictive distribution is absolutely continuous, it can be argued
that $\myS(x,Q)$ ought to depend only on the behavior of the
predictive density, $q$, in an infinitesimal neighborhood of the
observation that materializes, $x$. Any such scoring rule is said to
be local, with the logarithmic scoring rule,
%
\begin{equation} \label{eqLS}
\myS(x,Q) = - \ln q(x),
\end{equation}
being the most prominent example [\citet{Goo52}]. Another example is the
Hyv{\"a}rinen (\citeyear{Hyv05}) score,
%
\begin{eqnarray} \label{eqHS}
\myS(x,Q) & = & 2 \frac{q''(x)}{q(x)} - \biggl( \frac{q'(x)}{q(x)} \biggr)^{ 2}
\nonumber\\[-8pt]\\[-8pt]
& = & ((\ln q)'(x))^2 + 2 (\ln q)''(x),\nonumber
\end{eqnarray}
which is local of order 2, in the sense that it depends on the
predictive density only by its value, and the values of its first and
second derivative, at the observation. Similarly, the logarithmic
score can be considered to be local of order zero; in fact, it is the
only such score that is proper, up to equivalence [\citet{Ber79}].
The Hyv\"arinen score is also proper [\citet{DawLau05}], thus
raising the question for a~characterization of the local proper
scoring rules of order $k \leq2$.

In a~far-reaching recent paper, \citet{ParDawLau}
achieve a~characterization of the key local score functions of any
order $k \geq0$. They derive these scores from the Euler--Lagrange
equation of the calculus of variations, thereby obtaining natural
candidates for local proper scoring rules, the actual propriety of
which can be checked by additional criteria. We complement these
results---for more detailed comments, see Remark~\ref{re0}---by
developing an alternative approach, restricting ourselves to the
practically most relevant case of the local proper scoring rules of
order $k \leq2$. Our main contributions are the following: we build
on a~characterization of proper scoring rules via concave functionals
and their (super-)gradients, which yields the general form of the
second-order local proper scoring rules in a~natural tangent
construction; and we specify suitable classes of scoring rules and
predictive densities that allow for a~full-fledged, rigorous
characterization.

The remainder of the paper is organized as follows. Section
\ref{secrules} introduces the notions of propriety and locality in
full detail. Section~\ref{seccharacterization} presents our main
result, in that we characterize the class of the local scoring rules
of order~2 that are proper relative to a~comprehensive family of
Lebesgue densities, which includes many of the classical
location-scale families on the real line. In addition, we discuss the
relations to and distinctions from the work of \citet{ParDawLau}.
The proof of our main result is given in Section~\ref{secproof}.
Section~\ref{secancillary} provides supplements and examples, and a~data example on ensemble weather forecasts is given in Section
\ref{secdata}. Section~\ref{secdiscussion} closes with a~discussion
of open problems and hints at possible future developments and
applications.

\section{Local proper scoring rules} \label{secrules}

Initially, we consider predictive distributions on a~general sample space,
$\Omega$. Let $\cA$ be a~$\sigma$-algebra of subsets of $\Omega$, and
let $\cM$ be a~class of probability measures on $(\Omega,\cA)$.
A~function on $\Omega$ is $\cM$-\textit{quasi-integrable} if it is
measurable with respect to $\cA$ and quasi-integrable with respect to
all $Q \in\cM$ [\citet{Bau01}, page 64]. A~\textit{probabilistic forecast}
or a~\textit{predictive distribution} is any probability measure $Q \in
\cM$. A~\textit{scoring rule} is any extended real-valued function
$\myS\dvtx \Omega\times\cM\to\Real$ such that $\myS(\cdot,Q)$ is
$\cM$-quasi-integrable for all $Q \in\cM$. Hence, if the predictive
distribution is $Q$ and the event $\omega$ materializes, the
forecaster's loss is $\myS(\omega,Q)$. We define
\[
\myS(P,Q) = \int\myS(\omega,Q) \,d P(\omega)
\]
as the expected score under $P$ when the probabilistic forecast is
$Q$. This is a~well-defined extended real-valued quantity, because
$\myS(\cdot,Q)$ is quasi-integrable with respect to $P$.
%
\begin{definition} \label{defproper}
The scoring rule $\myS$ is \textit{proper} relative to $\cM$
if
\[
\myS(P,P) \leq\myS(P,Q) \qquad\mbox{for all } P, Q \in\cM.
\]
It is \textit{strictly proper} relative to $\cM$ if $\myS(P,P) \leq
\myS(P,Q)$ with equality if and only if $Q = P$.
\end{definition}

The term proper was coined by \citet{WinMur68}, while the
general idea can be traced to \citet{Bri50} and \citet{Goo52}.
\citet{Daw08} provides a~concise history of proper scoring rules, which
includes major contributions by the subjective school of probability
as well as meteorologists.

A~scoring rule can be thought of as \textit{local} if $\myS(\omega,Q)$
depends on the predictive distribution, $Q$, only through its behavior
in an infinitesimal neighborhood of the verifying observation,
$\omega$. Bernardo [(\citeyear{Ber79}), page 689] argued in this vein, noting that
``when assessing the worthiness of a~scientist's final conclusions,
only the probability he attaches to a~small interval containing the
true value should be taken into account.'' In the context of
predictive densities, the class $\cM$ is a~family of probability
measures that are absolutely continuous with respect to a~$\sigma$-finite measure $\mu$ on $(\Omega, \cA)$. We then identify a~probabilistic forecast $Q \in\cM$ with its $\mu$-density, $q$, which
we call a~\textit{predictive density} or a~\textit{density forecast}. The
classical example of a~local proper scoring rule is the aforementioned
logarithmic score, which can be interpreted as a~predictive
likelihood, and is strictly proper relative to any such class $\cM$.

Hereinafter, we restrict attention to the case in which the sample
space $\Omega$ is the real line, $\cA$ is the Borel $\sigma$-algebra,
$\mu$ is the Lebesgue measure, and $\cM$ corresponds to some class of
Borel probability measures that admit a~unique smooth Lebesgue
density, $q$. Accordingly, we will consider $\cM$ as a~class of
densities rather than measures, and we may write $\myS(\cdot,q)$. The
logarithmic score~(\ref{eqLS}) and the Hyv\"arinen score
(\ref{eqHS}) admit particularly simple analytic forms in terms of the
log-likelihood, $\ln q(x)$, and its derivatives, which are
fundamental objects of statistical inference. Therefore, we define
locality in terms of these quantities.
%
\begin{definition} \label{deflocal}
Let $k$ be a~nonnegative integer, and let $\cM$ be a~class of
probability densities with respect to the Lebesgue measure on $\real$
that are everywhere strictly positive and admit derivatives up to
order $k$. A~scoring rule $\myS$ for the class $\cM$ then is
\textit{local} of \textit{order} $k$ if there exists a~function $\mys\dvtx
\real^{2+k} \to\Real$, which we call a~\textit{scoring function}, such
that
\[
\myS(x,q)
= \mys\bigl( x, \ln q(x), \ldots, (\ln q)^{(k)}(x)\bigr)
\]
for every $q \in\cM$ and $x \in\real$.
\end{definition}

An alternative notion of locality, which allows the predictive
density,~$q$, to have zeroes, would take the arguments of the scoring
function as~$x, q(x), \ldots,\allowbreak q^{(k)}(x)$. However, in addition to
being natural and facilitating the technicalities, the assumption of
strict positivity avoids pathologies, as will be seen in Remark
\ref{re4} below.

As propriety can only be assessed relative to a~specified class of
predictive densities, we now introduce a~suitable family.
%
\begin{definition} \label{defP}
Let $\cP$ denote the class of all probability densities, $p$, with respect
to the Lebesgue measure on $\real$ that satisfy the following
conditions:
\begin{longlist}[(P4)]
\item[(P1)]
$p$ is strictly positive on $\real$;
\item[(P2)]
$p$ admits four continuous derivatives on $\real$;
\item[(P3)]
for every $m >0$ and $j = 0, 1, \ldots, 4$,
\[
{\lim_{x \rightarrow\pm\infty}} |x|^m p^{(j)}(x) = 0;
\]
\item[(P4)]
there exists a~constant $a~= a(p) > 0$ such that
\[
{ \lim_{x \rightarrow\pm\infty}} |x|^{-a} \frac{p^{(j)}(x)}{p(x)} = 0
\qquad\mbox{for } j = 1, \ldots, 4.
\]
\end{longlist}
\end{definition}

The class $\cP$ is quite broad and includes many well-known densities,
such as all normal and logistic densities, the corresponding skew
variants [\citet{autokey15}], and finite mixtures of these densities. In
particular, the class $\cP$ is convex, as implied by the following
result.
%
\begin{lemma} \label{leuint}
For every $k = 1, 2, \ldots$ there exists a~polynomial $M =
M(y_1,\allowbreak
\ldots, y_k)$ of degree $k$ such that for all $p, q \in\cP$ and
$\alpha\in[0,1]$, the density $r_\alpha= \alpha p +
(1-\alpha) q$ satisfies
\[
\bigl| (\ln r_\alpha)^{(k)}(x) \bigr| \leq M \biggl( \max \biggl\{
\frac{|p'(x)|}{p(x)}, \frac{|q'(x)|}{q(x)} \biggr\}, \ldots, \max
\biggl\{ \frac{|p^{(k)}(x)|}{p(x)}, \frac{|q^{(k)}(x)|}{q(x)} \biggr\}
\biggr)
\]
pointwise in $x \in\real$.
\end{lemma}
\begin{pf}
Let the polynomial $L(y_1,\ldots,y_k)$ of degree $k$ be such that the
$k$th logarithmic derivative of a~smooth function $g>0$ can be
written as
\[
(\ln g)^{(k)} = L \biggl( \frac{g'}{g}, \ldots, \frac{g^{(k)}}{g} \biggr) ,
\]
where here and in the following we suppress the argument $x \in
\real$. Define the polynomial $M$ as $L$ with all coefficients
replaced by their absolute values. Evidently then,
\[
\bigl|(\ln r_\alpha)^{(k)} \bigr| \leq
M \biggl( \frac{|r_\alpha'|}{r_\alpha}, \ldots, \frac{|r_\alpha
^{(k)}|}{r_\alpha} \biggr) ,
\]
and it suffices to show that
\[
\frac{|r_\alpha^{(j)}|}{r_\alpha} \leq\max
\biggl\{ \frac{|p^{(j)}|}{p}, \frac{|q^{(j)}|}{q} \biggr\}
\qquad\mbox{for } j = 1, \ldots, k.
\]
Consider the function $f(\alpha) = (\alpha c_1 + (1-\alpha)
c_0) / (\alpha d_1 + (1-\alpha) d_0)$, where~$c_0,\allowbreak c_1 \in
\real$ and $d_0, d_1 > 0$ are constants. Then
\[
|f(\alpha)| \leq\max\{|f(0)|, |f(1)|\} \qquad\mbox{for } \alpha\in[0,1],
\]
because $f'(\alpha) = (c_1 d_0 - c_0 d_1) / (\alpha d_1 +
(1-\alpha) d_0)^2$ does not change sign. The desired inequality
follows on setting $c_0 = q^{(j)}(x)$, $c_1 = p^{(j)}(x)$, $d_0 =
q(x)$ and $d_1 = p(x)$.
\end{pf}
%
%
\begin{corollary} \label{corP}
The class $\cP$ is convex.
\end{corollary}

In the following, we do not systematically distinguish a~scoring rule,
$\myS$, and the corresponding scoring function, $\mys$, both of which
will simply be referred to as scores.

\section{Characterization of the local proper scores of order 2}
\label{seccharacterization}

Along with the logarithmic and the Hyv\"arinen score, any convex
combination thereof is a~local proper score of order 2. However,
the class of the local proper scoring rules of order~2 on the real
line, $\real$, has a~much richer structure, and allows for a~characterization in terms of concave functionals.

\subsection{Main results} \label{secmain}

We first introduce classes of functions that satisfy suitable
polynomial growth conditions.
%
\begin{definition} \label{defRk}
Let $k$ be a~nonnegative integer. The class $\cR_k$ consists of all
functions $K\dvtx\real^{2+k} \to\real$ that admit continuous partial
derivatives up to order $2k$, and for which there exist finite
positive constants $C$ and $r$ such that, whenever $W$ stands for $K$
or any of its partial derivatives up to order~$2k$, then
\[
|W(x, y_0, \ldots, y_k)| \leq C
\{ (1+|x|) (1+|y_0|) \cdots(1+|y_k|) \}^r
\]
for all $(x, y_0, \ldots, y_k) \in\real^{2+k}$.
\end{definition}

Note that the growth conditions on the functions in the class $\cR_k$,
as well as the decay conditions on the densities in the class $\cP$ of
Definition~\ref{defP}, apply to each member individually. They are
not required to hold uniformly.

For a~function $K \in\cR_k$ and a~density $p \in\cP$, let
%
\begin{equation} \label{eqPhi}
\Phi_K( p) =
\int_\real K\bigl(x,\ln p(x),(\ln p)'(x),\ldots,(\ln p)^{(k)}(x)\bigr) p(x) \,dx.
\end{equation}
The integral exists and is finite by virtue of the growth and decay
conditions imposed on $K$ and $p$, respectively. Thus, any $K \in
\cR_k$ induces a~well-defined functional $\Phi_K\dvtx\cP\to\real$.
The role of the function $K$ in~(\ref{eqPhi}) resembles that of a~kernel in functional analysis. Hence, we will subsequently refer to
$K$ as a~\textit{kernel}, for ease of reference. The properties of such
kernels and the associated functionals play a~key role in our
subsequent characterization. In stating it, we use standard
abbreviations to denote the partial derivatives of a~function of the
form $g = g(x, y_0, \ldots, y_k)$; for example, we write $\pd_j g
= \pd g/ \pd y_j$ and $\pd^2_{xj} g = \pd^2 g / (\pd x \,\pd
y_j)$. The proof is given in Section~\ref{secproof}.

The subsequent two results are closely connected to the work of
\citet{ParDawLau}; see Remark~\ref{re0}.
%
\begin{theorem} \label{thcharacterization}
Let $\cP$ denote the class of probability densities introduced
in Definition~\ref{defP}.
\begin{longlist}[(c)]
\item[(a)]
Consider a~kernel $K$ of the form
%
\begin{equation} \label{eqK}
K(x, y_0, y_1) = c y_0 + K_0(x, y_1),
\end{equation}
where $c$ is a~real constant and $K_0$ is a~real function on
$\real^2$. If $K \in\cR_1$ and the functional $\Phi_K$ is
concave, the function $\mys\dvtx \real^4 \to\real$, defined by
%
\begin{equation} \label{eqs}
\mys(x,y_0,y_1,y_2) = c y_0 +
(1 - y_1 \pd_1 - \pd^2_{x1} - y_2 \pd^2_{11}) K_0(x,y_1),
\end{equation}
represents a~local score of order 2 that is proper relative
to $\cP$.
\item[(b)]
Conversely, if $\mys\in\cR_2$ represents a~local score of
order 2 that is proper relative to $\cP$, there exists a~kernel $K \in\cR_1$ of the form~(\ref{eqK}), where $c$ is a~real
constant and $K_0$ is a~real function on $\real^2$, such that the
functional $\Phi_K$ is concave and $\mys$ admits the representation
(\ref{eqs}).
\item[(c)]
The above statements remain valid with concave replaced by strictly
concave, and proper replaced by strictly proper.
\end{longlist}
\end{theorem}

The following sufficient condition for the functional $\Phi_K$ to be
concave will be proved in Section~\ref{secconcave}.
%
\begin{proposition} \label{propconcave}
Suppose that $K$ is a~kernel of the form~(\ref{eqK}) such that
\textup{(i)} $K \in\cR_1$, \textup{(ii)} $c \leq0$, and \textup{(iii)}
the map $y_1 \mapsto K_0(x, y_1)$ is concave for every $x \in\real$.
Then the functional $\Phi_K\dvtx\cP\to\real$ is concave. The statement
continues to hold if concave is replaced by strictly concave.
\end{proposition}

The criterion provides a~straightforward method of constructing local
proper scores of order 2 via the basic relationship~(\ref{eqs}).
For example, the kernel $K(x,y_0, y_1) = - y_0$ yields the scoring
function, $\mys(x, y_0, y_1, y_2) = - y_0$, that represents the
logarithmic score~(\ref{eqLS}). The associated functional
\[
\Phi( p) = \myS( p,p) = - \int_\real p(x) \ln p(x) \,dx
\]
is the Shannon entropy, and the associated divergence
\[
d_{\mathrm{KL}}( p, q) = \myS( p,q) - \myS( p,p)
= \int p(x) \ln\frac{p(x)}{q(x)} \,dx
\]
is the Kullback--Leibler divergence. Similarly, the kernel
$K(x,y_0,y_1) = - y_1^2$ yields the scoring function, $\mys(x, y_0,
y_1, y_2) = y_1^2 + 2y_2$, that represents the Hyv\"arinen score
(\ref{eqHS}). The associated functional and divergence
\[
-\int_\real\biggl( \frac{p'(x)}{p(x)} \biggr)^{ 2} p(x) \,dx
\quad\mbox{and}\quad
d_{\mathrm{FI}}( p, q) =
\int\biggl( \frac{p'(x)}{p(x)} - \frac{q'(x)}{q(x)} \biggr)^{ 2} p(x) \,dx
\]
are minus the Fisher information and the Fisher information distance
[\citet{Das08}, Definitions 2.5 and 2.6, pages 25 and 26], respectively.
For further examples, see Section~\ref{secexamples}.

\subsection{Remarks} \label{secremarks}

It has to be emphasized that the present work owes a~great deal to
interactions with Philip Dawid, Steffen Lauritzen and Matthew Parry,
which began with their kindly pointing out an error in our previous
work [\citet{EhmGne}].\vadjust{\goodbreak}
%
\begin{remark}[(Acknowledgment of priority)] \label{re0}
In the compact notation explained in Section~\ref{secproof}, any
second-order local proper scoring rule can be written as
%
\begin{equation} \label{eqPDL}
\mys = K - \biggl[ z_1 + \frac{d}{dx} \biggr] \pd_1 K.
\end{equation}
We learned about this representation in a~personal communication
[\citet{DawParLau}]. Detail on the relation of our work
to the paper by \citet{ParDawLau} is provided in the next remark.
\end{remark}
%
%
\begin{remark} \label{re01}
Employing an elegant approach based on operator algebra,
\citet{ParDawLau} investigate local proper scoring rules on a~general open
interval on the real line of any order $k \geq0$.
In a~\textit{tour de force}, they establish the existence of
key local score functions for any even order, and their nonexistence
for odd orders, in addition to studying their invariance under data
transformations. In the case $k = 2$ the general form~(39) of the
key local scoring rules in \citet{ParDawLau} is essentially
equivalent to ours, up to the parameterization in terms of densities
rather than log densities.

Despite the many parallels to the work of \citet{ParDawLau}, there
are important differences, including the basic approach and techniques
employed. A~key local score derives from the homogeneous
Euler--Lagrange equation, which characterizes the scores for which
every density $p$ is a~stationary point of the mapping $q \mapsto
\myS( p,q)$. Accordingly, Parry, Dawid and Lauritzen's (\citeyear{ParDawLau}) analysis is in
terms of differential calculus, which leads to separate discussions of
the boundary terms from partial integrations and of sufficient
conditions for (strict) propriety. The latter occur in Theorem 9.1 of
\citet{ParDawLau} in the form of concavity conditions on homogeneous
$q$-functions, which correspond to our kernels; Proposition~\ref{propconcave} states essentially the same result in the case $k =
2$.

In a~different ansatz, our work starts from the characterization of
proper scoring rules via concave functionals and their
\hbox{(super-)}gradients [\citet{HenBue71}, \citet{GneRaf07}]. This readily yields the basic form~(\ref{eqtconstrp})
of the second-order local proper scoring rules in a~natural tangent
construction, up to a~possibly nonlocal term. Only then we apply the
calculus of variations to show that the possibly nonlocal term
vanishes, which establishes the definite form~(\ref{eqs}). Control
of the boundary terms from partial integrations is vital, and is
achieved through our particular choice of the classes of scoring
functions and predictive densities. The explicit specification of the
classes $\cS$ and $\cD$, along with the tangent construction, allow us
to give a~rigorous, yet full-fledged and practically relevant
characterization of the second-order local proper scoring rules, hence
constitute the main original contributions of our work.
\end{remark}

We continue with comments relating to the choice of the class
$\cP$ and the complementary roles of the kernel $K$ as a~function and
a~functional, thereby touching on the generality of Theorem
\ref{thcharacterization} and Proposition~\ref{propconcave}.
%
\begin{remark} \label{re1}
There is a~slight asymmetry in Theorem~\ref{thcharacterization}, in
that under the conditions of the sufficiency part (a) the scoring
function is continuous only, whereas the necessity part (b) requires
it to be four times continuously differentiable. Other than this, the
theorem accomplishes a~full characterization of the local proper
scoring rules of order 2 relative to the class $\cP$ of Definition
\ref{defP}.
\end{remark}
%
%
\begin{remark} \label{re2}
Part (a) of Theorem~\ref{thcharacterization} expresses a~local proper
score of order~2, $\mys$, in terms of a~kernel, $K$, with suitable
properties. Similarly, part (b) admits a~constructive extension that
finds and expresses a~suitable kernel, $K$, in terms of a~local proper
score of order 2, $\mys$. See Section~\ref{secconstruction} for
the explicit construction and Example~\ref{expowerseries} for an
illustration.
\end{remark}
%
%
\begin{remark} \label{re3}
Theorem~\ref{thcharacterization} has been stated for the special
class $\cP$ of Definition~\ref{defP}. Propriety relative to such a~broad class is a~fairly demanding requirement, and from this
perspective, part (a) is a~strong result. In contrast, part (b) would
be stronger if propriety was required relative to a~subclass $\cP_0
\subset\cP$ only. On the other hand, $\cP_0$ must not be too narrow.
An inspection of Section~\ref{secproof} shows that part (b) remains
valid relative to any convex subclass $\cP_0 \subset\cP$ with the
following two additional properties:
\begin{longlist}[(P6)]
\item[(P5)]
if a~continuous function $f$ on $\real$ with at most polynomial growth
at $\pm\infty$ satisfies $\int_\real f(x) ( p(x)-q(x)) \,dx \geq
0$ for all $p, q \in\cP_0$, then $f$ is constant;
\item[(P6)]
the richness properties of Lemma~\ref{lerichness} hold for $\cP_0$.
\end{longlist}
The inequality in condition~(P5) can be replaced by equality,
making~(P5) a~variant of the classical property of completeness
of the family $\cP_0$. Property~(P5) is needed in Section~\ref{secvarcalc},
while property (P6) is required in Section~\ref{secreduction}. The
full class $\cP$ does satisfy these conditions.
\end{remark}
%
%
\begin{remark} \label{re4}
The sufficiency part of Theorem~\ref{thcharacterization} would be
stronger if the statement applied relative to larger classes $\cP_1
\supset\cP$. The following adaptation of an example of \citet{Hub74}
shows that any such extension may entail unexpected effects for strict
propriety, with undesirable consequences in applications. Suppose
that $\cP$ is augmented to a~convex class $\cP_1$ that includes the
densities
\[
p_\alpha(x) =
\cases{
\alpha g(x), &\quad if $x \geq0$, \cr
(1 - \alpha) g(-x), &\quad if $x < 0$,}
\]
where $\alpha\in(0,1)$ and $g(x) = x^5 e^{-x} / \Gamma(6)$
for $x \geq0$. The densities $p_\alpha$ satisfy all conditions for
the class $\cP$ except for property (P1), since $p_\alpha(x) = 0$ at
$x = 0$. As the logarithmic derivatives, $p'_\alpha(x)/p_\alpha(x)$,
do not depend on $\alpha$, the Fisher information of $p_\alpha$ does
not depend on $\alpha$ either, hence its negative is not strictly
concave as a~functional on $\cP_1$. Accordingly, the Fisher information
distance does not distinguish the densities $p_\alpha$, that is,
$d_{\mathrm{FI}}( p_\alpha, p_\beta) = 0$ for $\alpha, \beta\in
(0,1)$, and the Hyv\"arinen score~(\ref{eqHS}) fails to be strictly
proper relative to the augmented class $\cP_1$. In particular, strict
concavity of the function $K_0(x,y_1)$ of Proposition
\ref{propconcave} in $y_1$ does not imply strict concavity of the
associated functional, unless we restrict the class of densities under
consideration.
\end{remark}
%
%
\begin{remark} \label{re5}
$\!\!\!$By Proposition~\ref{propconcave}, concavity of a~kernel $K$ of the
form~(\ref{eqK}) in $y_1$ implies concavity of the associated
functional $\Phi_K$ on the class $\cP$. Conversely, what are the
consequences of concavity of the functional $\Phi_K$ on the kernel
$K$? The example of the logarithmic score~(\ref{eqLS}) demonstrates
that matters are not straightforward; here the functional $\Phi_K$ is
strictly concave, yet the kernel $K(x, y_0, y_1) = -y_0$ is not.

Now consider any kernel $K$ of the form~(\ref{eqK}) for which the
associated functional $\Phi_K$ is concave on $\cP$. Do we necessarily
have $c \leq0$ then? This is indeed true if $K_0(x, y_1) = - y_1^2$
represents the Hyv\"arinen score~(\ref{eqHS}). Then by propriety
\[
0 \leq\myS( p, q) - \myS(p,p) = d_{\mathrm{FI}}( p, q) - c d_{\mathrm{KL}}( p, q)
\]
for all $p, q \in\cP$, so that $c \leq0$ is necessary if the ratio
$r = d_{\mathrm{FI}}( p, q) / d_{\mathrm{KL}}( p, q)$ can attain
arbitrarily small values. However, $\cP$ contains all normal
densities, and if $p$ and $q$ are normal with mean zero and standard
deviations $\sigma$ and $\tau$, then
\[
d_{\mathrm{KL}}( p, q) = \frac12 \biggl[\frac{\sigma^2}{\tau^2} - 1 - \ln\frac
{\sigma^2}{\tau^2} \biggr]
\quad\mbox{and}\quad
d_{\mathrm{FI}}( p, q) = \frac{(1 - \sigma^2/\tau^2)^2}{\sigma^2},
\]
whence $r$ can attain any positive value. The argument clearly
depends on the class~$\cP$; it fails if $\cP$ is replaced by a~narrower class
$\cP_0$ for which the ratio~$r$ is bounded away from
zero. Such is in fact possible due to a~logarithmic Sobolev
inequality, which asserts that for certain classes $\cP_0 \subset\cP$
one has $d_{\mathrm{KL}}( p, q) \leq C d_{\mathrm{FI}}( p, q)$ for
$p, q \in\cP_0$ with a~constant $C$ that depends only on~$\cP_0$. A~corresponding reference is \citet{Vil09}: put $u = \sqrt{p/q}$ and
$d\nu(x) = q(x) \,dx$ in equation (21.3) and consider Remark 21.4.
\end{remark}

\section{\texorpdfstring{Proof of Theorem \protect\ref{thcharacterization}}{Proof of Theorem 3.2}} \label{secproof}

Our point of departure is Theorem 1 of \citet{GneRaf07},
which can be traced to \citet{HenBue71} and
characterizes proper scoring rules by means of the supergradients of
concave functionals on convex classes of probability measures. We
state it in the special case where that class corresponds to the set
$\cP$ of Lebesgue densities introduced in Definition~\ref{defP}.
\textit{Throughout this section propriety is understood as propriety
relative to~$\cP$.}
%
\begin{theorem} \label{thGR2007}
Let $\Phi$ be a~real-valued concave functional on $\cP$ with
supergradient $\Phi^*(\cdot,p)\dvtx\real\to\real$ at $p \in\cP$,
that is,
\[
\Phi( q) - \Phi(p) - \int_\real\Phi^*(x,p) \bigl( q(x) - p(x) \bigr) \,dx \leq0
\qquad\mbox{for } p, q \in\cP.
\]
Then the scoring rule
%
\begin{equation} \label{eqspropform}
\myS(\cdot,p) = \Phi^*(\cdot,p) - \int_\real\Phi^*(x,p) p(x) \,dx + \Phi(p)
\end{equation}
is proper, and
\[
\myS(p,p) = \int_\real\myS(x,p) p(x) \,dx = \Phi(p)
\qquad\mbox{for } p \in\cP.
\]
Conversely, if $\myS$ is proper, then $\Phi(p) = \myS(p,p)$ is a~concave functional on $\cP$ with supergradient $\Phi^*(\cdot,p) =
\myS(\cdot,p)$ at $p \in\cP$, whence $\myS$ is of the form
(\ref{eqspropform}). Furthermore, the above continues to hold with
concave replaced by strictly concave, and proper replaced by strictly
proper.
\end{theorem}

For sufficiently regular local proper scoring rules we can compute
gradients of the corresponding functionals. Specifically, a~function
$G(\cdot,p)\dvtx\real\to\real$ is a~\textit{weak gradient}, or simply a~\textit{gradient}, of the functional $\Phi$ at $p \in\cP$ if for every
$q \in\cP$
%
\begin{equation} \label{deftangent}
\frac{d}{dt} [ \Phi( q_t) ] \bigg|_{t=0}
= \int_\real G(x,p) \bigl( q(x) - p(x) \bigr) \,dx,
\end{equation}
where
%
\begin{equation} \label{eqmixture}
q_t = (1-t) p + t q \qquad\mbox{for } t \in[0,1].
\end{equation}
Any (super-)gradient is defined only modulo an arbitrary additive
constant that may depend on $p$, which does not affect the
construction~(\ref{eqspropform}).

Theorem~\ref{thGR2007} along with such tangent calculations gives us
a~construction method for local proper scores that readily elucidates
their particular form. We refer to this approach as the \textit{tangent
construction} and give details in the following section, before
completing the proof of Theorem~\ref{thcharacterization} in a~series
of subsequent steps.

\subsection{Tangent construction of proper scores} \label{sectangent}

In what follows we use compactified notation whenever possible. As
noted, we do not systematically distinguish scoring rules, $\myS$, and
the corresponding scoring functions, $\mys$, both of which are
referred to as scores. Log-likelihoods and their derivatives are
denoted by $z_0(x,p) = \ln p(x)$ and
\[
z_j(x,p) = (\ln p)^{(j)}(x) \qquad\mbox{for } j = 1, 2, \ldots
\]
or simply $z_0$ and $z_j$ if the density $p$ is fixed. Clearly then,
\[
z_j' = z_{j+1} = z_0^{(j+1)} \qquad\mbox{for } j = 0, 1, 2, \ldots,
\]
where the prime denotes differentiation with respect to $x$. We
usually suppress the differential, $dx$, in integrals over $x \in
\real$, and in the corresponding integrands we omit all or part of the
arguments whenever these are clear from the context. For example,
given $K \in\cR_k$ and $p, q \in\cP$ we may abbreviate\looseness=-1
\[
\int_\real K\bigl(x,\ln q(x),\ldots,(\ln q)^{(k)}(x)\bigr) p(x) \,dx
\]\looseness=0
as
\[
\int K p\qquad (K = K_q),
\]
where, evidently, $K_q = K_q(x) = K(x,\ln q(x),\ldots,(\ln q)^{(k)}(x))$.

We now develop the tangent construction. The first step consists in
calculating the gradients of (not necessarily concave) functionals of
kernel type.\vspace*{-3pt}

\begin{lemma} \label{lemgrad2}
Let $K \in\cR_2$. Then a~gradient $G$ of the associated
functional $\Phi_K\dvtx\cP\to\real$ exists at any $p \in\cP$
and is given uniquely by
%
\begin{equation} \label{eqgrad2}
G = K + \pd_0 K - \frac{1}{p} \,\frac{d}{dx} [ p \pd_1 K ]
+ \frac{1}{p} \,\frac{d^2}{dx^2} [ p \pd_2 K ]\qquad (G = G_p, K =
K_p)\hspace*{-23pt}
\end{equation}
up to an arbitrary additive constant that may depend on $p$.\vspace*{-3pt}
\end{lemma}

Recall that according to our notational conventions,
(\ref{eqgrad2}) means that the relation holds whenever the
functions $G$, $K$ and $\pd_j K$ are evaluated at arguments
\[
(x,z_0,z_1,z_2) = (x,\ln p(x),(\ln p)'(x),(\ln p)''(x)),
\]
where $p \in\cP$ and $ x \in\real$.\vspace*{-3pt}
\begin{pf*}{Proof of Lemma~\ref{lemgrad2}}
Let $p \in\cP$ be fixed. In calculating a~gradient of $\Phi=
\Phi_K$ at $p$ we initially ignore all technicalities, that is, we
assume that integrals are well defined and finite, that the order of
integration and differentiation can be interchanged, and that boundary
terms in partial integrations vanish. Then
%
\begin{equation} \label{eqgcalc1}
\frac{d}{dt} [ \Phi( q_t) ]
= \int\frac{d}{dt} [ K_t q_t ]
= \int K_t (q-p) + \int\biggl[ \frac{d}{dt} K_t \biggr] q_t,
\end{equation}
where $q_t$ denotes the mixture density~(\ref{eqmixture}) and
\[
K_t = K_{q_t} = K( x, \ln q_t(x), (\ln q_t)'(x), (\ln q_t)''(x)).
\]
Since $\frac{d}{dt} \ln q_t = ( q - p) / q_t$, the mixed
derivative with respect to $t$ and $x$ of order $j$ is given by
\[
\frac{d}{\,dt} \bigl[(\ln q_t)^{(j)} \bigr] = \biggl( \frac{q-p}{q_t} \biggr)^{(j)}.\vadjust{\goodbreak}
\]
The second term on the right-hand side of~(\ref{eqgcalc1}) can
then be computed using partial integration, in that
%
\begin{eqnarray} \label{eqgcalc2}
&& \int\biggl[ \frac{d}{dt} K_t \biggr] q_t \nonumber\\
&&\qquad = \int\biggl[ (\pd_0 K_t) \biggl( \frac{q-p}{q_t} \biggr)
+ (\pd_1 K_t) \biggl( \frac{q-p}{q_t} \biggr)'
+ (\pd_2 K_t) \biggl( \frac{q-p}{q_t} \biggr)'' \biggr] q_t \nonumber\\
&&\qquad = \int(\pd_0 K_t) ( q - p)
- \int\biggl( \frac{d}{dx} [ q_t \pd_1 K_t] \biggr) \biggl( \frac{q-p}{q_t} \biggr)
\\
&&\qquad\quad{}
+ \int\biggl( \frac{d^2}{dx^2} [ q_t \pd_2 K_t] \biggr) \biggl( \frac{q-p}{q_t} \biggr)
\nonumber\\
&&\qquad = \int\biggl[ \pd_0 K_t
- \frac{1}{q_t} \,\frac{d}{dx} [ q_t \pd_1 K_t ]
+ \frac{1}{q_t} \,\frac{d^2}{dx^2} [ q_t \pd_2 K_t ] \biggr] ( q - p).
\nonumber
\end{eqnarray}
Evaluating at $t = 0$, and noting that $q_0 = p$ and $K_0 = K_p = K$,
(\ref{eqgcalc1}) and~(\ref{eqgcalc2}) yield
\[
\frac{d}{dt} [ \Phi( q_t) ] \bigg|_{t=0}
= \int\biggl[ K + \pd_0 K - \frac{1}{p} \,\frac{d}{dx} [ p \pd_1 K ]
+ \frac{1}{p} \,\frac{d^2}{dx^2} [ p \pd_2 K ] \biggr] ( q - p),
\]
showing that $G$ from~(\ref{eqgrad2}) is indeed a~gradient of $\Phi$ at $p$.

It remains to settle the technicalities. Generally, if a~family $\{
h(x,t)\dvtx x \in\real,\allowbreak t \in[0,1] \}$ is such that $h(x,t)$ is
integrable with respect to $x$ for every $t$, and the family $\{ \pd_t
h(x,t)\dvtx x \in\real, t \in[0,1] \}$ of partial derivatives is
uniformly integrable and continuous in $t$ for every $x$, then $H(t) =
\int h(x,t) \,dx$ is differentiable with
\[
\frac{dH}{dt}(0) = \int\pd_t h(x,0) \,dx.
\]
Here we consider~(\ref{eqgcalc1}) and identify $h(\cdot,t) = K_t q_t$,
so that
%
\begin{eqnarray} \label{eqhdt}
\pd_t h(\cdot, t) &=& (K_t + \pd_0 K_t) ( q - p) \nonumber\\[-8pt]\\[-8pt]
&&{}
+ \biggl[ (\pd_1 K_t) \biggl( \frac{q-p}{q_t} \biggr)'
+ (\pd_2 K_t) \biggl( \frac{q-p}{q_t} \biggr)'' \biggr] q_t. \nonumber
\end{eqnarray}
Now $\pd_t h(\cdot,t)$ is continuous in $t$, because $K$ and its
partial derivatives are continuous, and their arguments depend
continuously on $t$. Concerning uniform integrability, each of the
terms $K_t$, $\pd_0 K_t$, $\pd_1 K_t$, $\pd_2 K_t$ grows at most
polynomially as $x \to\pm\infty$. This is because by Lemma
\ref{leuint} and property (P4) of the class $\cP$ the arguments of
the terms grow at most polynomially; as $K \in\cR_2$, the same is
true for the functions themselves. Furthermore, by property (P3) and
the above, the terms $q - p$,
\[
\biggl( \frac{q-p}{q_t} \biggr)' q_t = \biggl( \frac{q'-p'}{q_t} -
\frac{q-p}{q_t} \frac{q_t'}{q_t} \biggr) q_t = q' - p' - ( q - p) (\ln
q_t)'
\]
and
\begin{eqnarray*}
\biggl( \frac{q-p}{q_t} \biggr)'' q_t & = &
q'' - p'' - 2 ( q' - p') (\ln q_t)' - ( q - p) (\ln q_t)'' \\
&&{} + ( q - p) ((\ln q_t)' )^2,
\end{eqnarray*}
decay faster than the
reciprocal of any polynomial as $x \rightarrow\pm\infty$.
Therefore, the corresponding products in~(\ref{eqhdt}) with the terms
involving $K$ decay faster than the
reciprocal of any polynomial as well. By Lemma
\ref{leuint}, this property holds uniformly in $t \in[0,1]$. Thus,
the family~(\ref{eqhdt}) is uniformly integrable, and we may
interchange the order of the integration and differentiation. Similar
growth and decay considerations show that the boundary terms in the
partial integrations in~(\ref{eqgcalc2}) vanish.

Finally, uniqueness follows from the property (P5) satisfied by the
class~$\cP$ (cf. Remark~\ref{re3}) and the at most polynomial growth
of $G$ as $|x| \to\infty$.
\end{pf*}

For use later on, we also state a~version of Lemma~\ref{lemgrad2}, in
which $K \in\cR_1$ so that $\pd_2 K$ vanishes. The proof is
analogous.
%
\begin{lemma} \label{legrad1}
Suppose that the kernel $K$ depends on arguments $x, z_0$ and~$z_1$ only and belongs to $\cR_1$. Then a~gradient $G$ of the
associated functional $\Phi_K\dvtx\cP\to\real$ exists at any $p
\in\cP$, and is given uniquely by
%
\begin{equation} \label{eqgrad1}
G = K + \pd_0 K - \frac{1}{p} \,\frac{d}{dx} [ p \pd_1 K ]\qquad
(G = G_p, K = K_p)
\end{equation}
up to an arbitrary additive constant that may depend on $p$.
\end{lemma}

Hereafter we will ignore the irrelevant additive constant and refer to
the expression in~(\ref{eqgrad2}) and~(\ref{eqgrad1}), respectively,
as \textit{the tangent} of $\Phi$ at $p$. A~common form of the tangent
valid for both $k = 1$ and $k = 2$ is
%
\begin{equation} \label{eqcommon}
G = K + \pd_0 K + L_0 K,
\end{equation}
where the differential operator $L_0$ is formally defined via the
infinite sum
%
\begin{equation} \label{eqL0}
L_0K = \sum_{j=1}^\infty (-1)^j \frac{1}{p} \,\frac{d^j}{dx^j} [ p \pd_j K],
\end{equation}
and $L_0$ and $K$ depend tacitly on $p$. If $K \in\cR_k$, all but
the first $k$ terms in the sum vanish, and so the definition makes
good sense. In terms of the operator $L$ in equation (19) of
\citet{ParDawLau} we have $L = - p (\pd_0 + L_0)$.

For the second step of the tangent construction let again $\Phi=
\Phi_K$ be a~kernel type functional associated with some kernel $K$ in
$\cR_1$ or $\cR_2$. If~$\Phi$ is concave, then the tangent $G$ of
$\Phi$ at $p \in\cP$ is easily seen to be also a~supergradient, and by
Theorem~\ref{thGR2007} a~proper score is obtained by setting
%
\begin{eqnarray}\label{eqtconstr}
\mys & = & G - \int G p + \Phi(p)\qquad
(\mys= \mys_p, G = G_p) \nonumber\\[-8pt]\\[-8pt]
& = & K + \pd_0 K + L_0 K - \int (\pd_0 K) p\qquad
(K = K_p, L_0 = L_{0,p}).\nonumber
\end{eqnarray}
As for the step leading to~(\ref{eqtconstr}), note that in view of
(\ref{eqPhi}) and~(\ref{eqcommon}) we have
\[
\Phi(p) - \int G p
= \int K p -
\biggl( \int K p + \int (\pd_0 K) p + \int (L_0 K) p \biggr)
= - \int (\pd_0 K) p
\]
on using the fact that $\int(L_0K) p = 0$. This latter equality
holds because the integrand is a~total derivative, the primitive of
which vanishes as $x \rightarrow\pm\infty$, due to the growth and
decay properties of the functions in $\cR_k$ and densities in $\cP$.
We will refer to this trivial observation as the \textit{vanishing
argument}. Since it will be used several times we state it as a~lemma, despite its simplicity.
%
\begin{lemma}[(Vanishing argument)] \label{levanarg}
Let $p \in\cP$, and let $W$ be a~real, differentiable function
such that the function $g = p^{-1} \,\frac{d}{dx} [p W ]$ is
$p$-integrable, $\int|g| p < \infty$, and $\lim_{x \to\pm\infty}
p(x)W(x) = 0$. Then $\int g p = 0$.
\end{lemma}

This concludes the tangent construction of a~proper score $\mys$ from
a~concave functional $\Phi= \Phi_K$ where $K \in\cR_k$ ($k = 1, 2$).
We summarize the foregoing discussion.
%
\begin{proposition}[(Tangent construction)] \label{proptconstr}
Suppose that $K \in\cR_k$ where $k=1$ or $k=2$, and that the
associated functional $\Phi_K$ is concave. Then
%
\begin{equation} \label{eqtconstrp}\quad
\mys = K + L_0 K + \pd_0 K - \int(\pd_0 K) p\qquad
(\mys= \mys_p, K = K_p, L_0 = L_{0,p})
\end{equation}
is a~proper score relative to $\cP$. It is local of order $2k$
if $\pd_0 K$ is constant in~$x$ for every $p \in\cP$, or
if $\int(\pd_0 K) p$ does not depend on $p$.
\end{proposition}
\begin{pf}
The first claim has already been proved. Locality under the stated
conditions is obvious: if $\pd_0 K = \pd_0 K(x,\ln p(x),(\ln p)'(x))$
(for $k = 1$, say) does not depend on $x$, it equals its expectation,
$\int(\pd_0 K) p $. Finally, an explicit evaluation of the total
differential(s) in the term $L_0 K$ yields partial derivatives of
order $\leq2k$ only, thereby proving the order $2k$ claim.
\end{pf}

\subsection{Variational calculus} \label{secvarcalc}

A~score $\mys$ is proper if the functional $\cP\ni q \mapsto\myS(
p, q)$ achieves its minimum at $q = p$, \textit{for every} $p
\in\cP$. This circumstance allows a~variational characterization
of---in fact, a~necessary condition for---the local proper scores.\vadjust{\goodbreak}
%
\begin{lemma} \label{levarchar}
Suppose that $\mys\in\cR_2$ is a~local proper score relative to $\cP$.
Then for every $p \in\cP$ one has
%
\begin{eqnarray} \label{eqeeqn}
\pd_0 \mys+ L_0 \mys
= \pd_0 \mys- \frac{1}{p} \,\frac{d}{dx} [ p \pd_1 \mys]
+ \frac{1}{p} \,\frac{d^2}{dx^2} [ p \pd_2 \mys]
= c_p \nonumber\\[-8pt]\\[-8pt]
&&\eqntext{(s = s_p, L_0 = L_{0,p})}
\end{eqnarray}
on $\real$, where
%
\begin{equation} \label{eqcqexpl}
c_p = \int(\pd_0 \mys) p\qquad (s = s_p).
\end{equation}
\end{lemma}
\begin{pf}
Fix $p \in\cP$ and consider convex combinations of the form $q_t =
(1-t)p + tq$ where $q \in\cP$ and $t \in(0,1)$. As the score is
proper, we have $(\myS( p,q_t) - \myS( p,p))/t \geq0$ for
every $t$. Let us compute the limit as $t \rightarrow0$. Putting
$\mys_t = \mys_{q_t}$, we have at first
\[
t^{-1} \bigl(\myS( p,q_t) - \myS( q_t,q_t) \bigr)
= t^{-1} \int\mys_t ( p - q_t)
= - \int\mys_t ( q - p).
\]
Arguing in the same way as in the proof of Lemma~\ref{lemgrad2}, we
find that the integrand is uniformly integrable and continuous in $t$,
so the limit exists and equals $-\int\mys ( q - p)$ where
$\mys= \mys_0 = \mys_p$. Thus writing
\[
\myS( p,q_t) - \myS( p,p)
= \myS( p,q_t) - \myS( q_t,q_t) + \myS( q_t,q_t) - \myS(p,p)
\]
and using Lemma~\ref{lemgrad2} and~(\ref{eqL0}), we get
\begin{eqnarray*}
&& \lim_{t \rightarrow0} t^{-1} \bigl(\myS( p,q_t) - \myS( p,p)\bigr) \\
&&\qquad = - \int \mys ( q-p) + \int
\biggl( \mys+ \pd_0 \mys- \frac{1}{p} \,\frac{d}{dx}
[ p \pd_1 \mys] + \frac{1}{p} \,\frac{d^2}{dx^2}
[ p \pd_2 \mys] \biggr) ( q-p) \\
&&\qquad = \int ( \pd_0 \mys+ L_0 \mys) ( q-p) .
\end{eqnarray*}
It follows that
%
\begin{equation} \label{eqvarcneg}
\int ( \pd_0 \mys+ L_0 \mys) ( q-p) \geq0
\end{equation}
for every $q \in\cP$. We proceed to show that this is possible only if
$\pd_0 \mys+ L_0 \mys$ equals some constant $c_p$ almost everywhere,
hence everywhere by continuity. To this end, let $f = \pd_0 \mys+
L_0 \mys$ and $g = f - \int f p$. Then $\int g q \geq0$
for every $q \in\cP$. Suppose $g$ were not constant. Since $\int g
p = 0$, the Lebesgue measure of the (open) set $\{g < 0\}$ is
strictly positive. Thus, there exists a~probability density $q_1 \in
C^\infty$ with compact support such that $\int g q_1 < 0$. Then
$q = \frac{1}{2} ( q_1 + p) \in\cP$ and $\int g q =
\frac{1}{2} \int g q_1 < 0$, in contradiction to
(\ref{eqvarcneg}). Finally, the constant $c_p$ is easily identified
by integrating~(\ref{eqeeqn}) against $p$ and noting that $\int(L_0
\mys) p = 0$, by the vanishing argument.\vadjust{\goodbreak}
\end{pf}

Equation~(\ref{eqeeqn}) essentially is the Euler equation of the
calculus of variations [\citet{GelFom63}, pages 40--42] and
corresponds to equation~(24) of \citet{ParDawLau}. Its slightly
different form here results from the fact that in our case the
integrand of the functional to be optimized is of the form $F(x, \ln
y, (\ln y)', (\ln y)'')$ rather than of the common form $F(x, y, y',
y'')$.

As a~first application of the Euler equation we show that local proper
scores are fixed points of the tangent construction. To this end, let
$\mys\in\cR_2$ be a~local proper score of order 2. By Theorem
\ref{thGR2007} and Lemma~\ref{lemgrad2}, the functional $\Phi_\mys\dvtx
\cP\to\real$ associated with the kernel $\mys$ is concave. The
tangent construction then gives
%
\begin{equation} \label{eqss}
\widetilde\mys= \mys+ L_0 \mys+ \pd_0 \mys - \int(\pd_0 \mys) p
\end{equation}
on substituting $\mys$ for $K$ in Proposition~\ref{proptconstr}.
Initially, this is another proper score, possibly of higher order, and
possibly nonlocal. However, by Lemma~\ref{levarchar} the right-hand
side of~(\ref{eqss}) reduces to $\mys$, whence in fact $\widetilde
\mys= \mys$.
%
\begin{proposition} \label{propinvar}
For a~local proper score $\mys\in\cR_2$ the tangent construction
based on the (concave) functional $\Phi_\mys$ leads back to $\mys$.
That is, any local proper score of order 2 is a~fixed point of the
tangent construction.
\end{proposition}

\subsection{Construction of a~$z_2$-independent kernel} \label{secconstruction}

The vanishing argument of Lemma~\ref{levanarg} enables us to modify a~given kernel without changing the associated functional. This
strategy is utilized in the following explicit construction of a~$z_2$-independent kernel from a~given local proper score. It is
analogous to the idea of gauge choice developed in Sections 7.3 and
7.4 of \citet{ParDawLau}. Again, it is tacitly assumed that the
quantities $z_j$ refer to a~fixed density $p \in\cP$, that is, $z_j =
z_j(x,p)$. As before, we frequently suppress these quantities when
they serve merely as arguments.
%
\begin{proposition} \label{propR1kernel}
Given the local proper score $\mys\in\cR_2$, let the kernel $K$
be defined as
%
\begin{equation} \label{eqKdef}
K = \mys- \frac{1}{p} \,\frac{d}{dx} [ p V]
= \mys- \biggl[ z_1 + \frac{d}{dx} \biggr] V,
\end{equation}
where
%
\begin{equation} \label{eqVdef}
V = \int_0^{z_1} \pd_2 \mys(x,z_0,t,z_2) \,dt.
\end{equation}
Then $K \in\cR_1$ and $\Phi_K = \Phi_\mys$. In particular, the
score $\mys$ can be reconstructed from the kernel $K$ via the
tangent construction.
\end{proposition}
\begin{pf}
The kernel $K$ inherits the polynomial growth properties from~$\mys$,
and it is twice continuously\vadjust{\goodbreak} differentiable since $\mys\in
C^4$. In particular, $K$ is well defined. An application of the
vanishing argument to the term $\frac{1}{p} \,\frac{d}{dx} [p
V]$ shows that $K$ and $\mys$ give rise to the same
functional. Thus $\Phi_K = \Phi_\mys$, and the last claim follows
from Propositions~\ref{proptconstr} and~\ref{propinvar}. Therefore,
to complete the proof it remains to show that the kernel $K$ from
(\ref{eqKdef}) does not depend on $z_2$, that is, we need to show
that $\pd_2 K = 0$.

A~comparison of the two differential operators
\begin{eqnarray*}
\pd_2 \,\frac{d}{dx} & = & \pd_2 ( \pd_x + z_1 \pd_0 + z_2 \pd_1 + z_3 \pd
_2 ) \\
& = & \pd_{x2}^2 + z_1 \pd_{02}^2 + \pd_1 + z_2 \pd_{12}^2 + z_3 \pd_{22}^2
\end{eqnarray*}
and
\[
\frac{d}{dx} \pd_2 = \pd_{x2}^2 + z_1 \pd_{02}^2 + z_2 \pd_{12}^2 + z_3
\pd_{22}^2
\]
yields the commutation relation
%
\begin{equation} \label{eqcommrel2}
\pd_2 \,\frac{d}{dx} = \frac{d}{dx} \pd_2 + \pd_1.
\end{equation}
Therefore, using $\pd_1 V = \pd_2 \mys$ [see~(\ref{eqVdef})] we get
\[
\pd_2 K
= \pd_2 \mys- z_1 \pd_2 V - \pd_1 V - \frac{d}{dx}\pd_2 V
= - \biggl[z_1 + \frac{d}{dx}\biggr]\pd_2 V,
\]
where
%
\begin{equation} \label{eqd2V}
\pd_2 V = \int_0^{z_1} \pd_{22}^2 \mys(x,z_0,t,z_2) \,dt.
\end{equation}
Thus if
%
\begin{equation} \label{eqvanish}
\pd_{22}^2 \mys(x,z_0(x,p),t,z_2(x,p)) = 0 \qquad\mbox{for all $x \in
\real, p \in\cP, |t| \le|z_1(x,p)|$,}\hspace*{-35pt}
\end{equation}
then $\pd_2 V=0$ and hence $\pd_2 K=0$, that is, $K \in\cR_1$, as
claimed. The somewhat lengthy proof of~(\ref{eqvanish}) is given in
the next subsection.

\subsection{\texorpdfstring{Proof of the independence condition (\protect\ref{eqvanish})}{Proof of the independence condition (27)}}
\label{secreduction}

The proof primarily rests upon the Euler equation~(\ref{eqeeqn}). An
evaluation of the total derivatives in~(\ref{eqeeqn}) shows that the
Euler equation can be written in the form
%
\begin{equation} \label{eqexpl1}
z_0^{(4)} \cdot a(x,z_0,z_0',z_0'',z_0''' ) -
b(x,z_0,z_0',z_0'',z_0''' ) = c_p\qquad (z_0 = \ln p)\hspace*{-35pt}
\end{equation}
with
\[
a~= \pd^2_{22} \mys \qquad\mbox{[so that in fact $a~= a(x,z_0,z_0',z_0'')$]}
\]
and a~function $b$, which depends (only) on the scoring function
$\mys$ and its partial derivatives up to order 3, other than $x$,
$z_0$ and the logarithmic derivatives~$z_1, z_2$ and $z_3 = z_2'$.
Therefore, and because $\mys$ is of smoothness class $C^4$, the
function $b$ is continuously differentiable. The same holds for $a~=
\pd^2_{22}\mys$, of course, so that $a, b \in C^1$.\vadjust{\goodbreak}

A~step critical to the remainder of the proof consists in showing that
the constant $c_p$ is independent of $p$. We state this below as
Proposition~\ref{propcconst}; its proof hinges on an argument due to
\citet{ParDawLau}. The ensuing fact that one and the same equation,
(\ref{eqexpl1}) with $c_p = c$, holds for every $p \in\cP$ is then
utilized to complete the proof of~(\ref{eqvanish}). For each of
these steps it is important that the class $\cP$ be sufficiently rich.
Let
\[
Z(x,q,k) = (z_0(x,q),z_1(x,q),\ldots,z_k(x,q))
\]
for $k = 0, 1, \ldots$ with $z_j(x,q) = (\ln q)^{(j)}(x)$ as above.
%
\begin{lemma}[($\cP$-richness)] \label{lerichness}
Let $k \in\{0, 1, \ldots, 4\}$. \textup{(a)} For every $x \in\real$
and $y \in\real^{k+1}$ there exists $q \in\cP$ such that
$Z(x,q,k) = y$. \textup{(b)} For every pair $q_1, q_2 \in\cP$ there
exist $q \in\cP$ and $x \in\real$ such that $Z(x,q,k) =
Z(x,q_1,k)$ and $q(u) = q_2(u)$ for $u$ outside some neighborhood
of $x$.
\end{lemma}
\begin{pf}
This is fairly obvious from the definition of $\cP$. For
completeness, we include a~proof. As concerns part (a), let $x \in
\real$ and $y \in\real^{k+1}$ be fixed. There is some $p \in\cP$
such that $z_0(x,p) = y_0$. We will construct $q$ as a~perturbation
$q = p (1+\psi)$ of $p$ such that $\psi(x)=0$. Certainly $q \in\cP$
if $\psi$ has compact support and is such that $1 + \psi> 0$, $\int
\psi p = 0$, and $\psi\in C^4$. It suffices to show that it is
possible to prescribe arbitrary values for the first $k$ derivatives
of $\psi$ at $x$, subject to those conditions. To that end, let $\psi
= Q M$ in the sense that $\psi(u) = Q(u) M(u)$ for $u \in\real$,
where
\[
Q(u) = \sum_{j=1}^r a_j (u-x)^j
\]
is a~polynomial vanishing at $x$ and $0 \le M \in C^4$ is a~mollifier
type function with (small) compact support $S$ such that $M(x) = 1$
and $M^{(j)}(x) = 0$ for $j = 1, \ldots, k$. Then $\psi^{(j)}(x) =
Q^{(j)}(x)$ for $j = 1, \ldots, k$, thereby confirming that arbitrary
values can indeed be prescribed if $r \geq k$. By increasing $r$ if
necessary, one can further assume that $Q$ attains both positive and
negative values on $S$.

Let any such $Q$ be fixed. We show that the conditions $Q M > -1$ and
$\int Q M p = 0$ can be satisfied, too. In fact, since $Q(x)=0$ one
can modify $M$ such that $Q M > -1$ everywhere without affecting its
local behavior at $x$. Since\vspace*{1pt} $\int_{\{Q < 0\}} Q M p < 0$ there is
$\delta>0$ such that the interval $J = [ x-\delta, x+\delta]$ is
contained in the interior of $S$ and $\int_J Q M p + \int_{J^c \cap
\{Q < 0\}} Q M p < 0$. Finally, on the set $S \cap J^c \cap\{Q >
0\}$ one can modify $M$ such that
\[
\int Q M p
= \int_J Q M p + \int_{J^c \cap\{Q < 0\}} Q M p + \int_{J^c \cap\{Q >
0\}} Q M p
= 0
\]
without affecting the condition $Q M > -1$. This concludes the proof
of (a). For part (b), note that because $q_1$, $q_2$ are continuous
probability densities,\vadjust{\goodbreak} there is $x \in\real$ such that $q_1(x) =
q_2(x)$. The above construction then yields a~local perturbation $q
\in\cP$ of $q_2$ satisfying $Z(x,q,k) = Z(x,q_1,k)$.
\end{pf}
%
%
\begin{proposition} \label{propcconst}
The constant $c_p$ in~(\ref{eqexpl1}) does not depend on the
density~$p$: there is a~constant $c \in\real$ such that $c_p =
c$ for every $p \in\cP$.
\end{proposition}
\begin{pf}
We use an argument in Section 4 (around Condition 4.1) of
\citet{ParDawLau}. Equation~(\ref{eqexpl1}) [resp.,
(\ref{eqeeqn})] can be condensed to a~statement of the form
%
\begin{equation} \label{eqexpl2}
F(x,Z(x,p,4)) = c_p\qquad (x \in\real, p \in\cP),
\end{equation}
where the function $F$ is determined by the score $\mys$ alone. Suppose
that $c_p$ is not independent of $p$. Then there are $q_1, q_2 \in\cP$
such that $c_{q_1} \neq c_{q_2}$. By Lem\-ma~\ref{lerichness}(b) there
exist $p \in\cP$ and $x_1 \neq x_2 \in \real$ such that $Z(x_1,p,4) =
Z(x_1,q_1,4)$ and $Z(x_2,p,4) = Z(x_2,q_2,4)$. By~(\ref{eqexpl2}) it
follows that for both $j = 1$ and $j = 2$
\[
c_{q_j} = F(x_j,Z(x_j,q_j,4)) = F(x_j,Z(x_j,p,4)) = c_p.
\]
The contradiction implies that $c_p$ is indeed independent of $p$.
\end{pf}

The following lemma is an easy consequence of the uniqueness theorem
for higher-order differential equations.
%
\begin{lemma}[(Reduction principle)] \label{lereduction}
Let $k \in\{0,1,2,3\}$, and let $a$ and $b$ be $C^1$
functions of arguments $x, y_0, \ldots, y_k$. Suppose that the
function $z_0 = \ln p(x)$, $x \in\real$ is, \textit{for every $p \in
\cP$}, a~solution of the differential equation
%
\begin{equation} \label{eqgenode}
z^{(k+1)} \cdot a\bigl(x,z,\ldots,z^{(k)}\bigr) = b\bigl(x,z,\ldots,z^{(k)}\bigr).
\end{equation}
Then $a(x,Z(x,p,k)) = 0$ for every $x \in\real, p \in\cP$.
\end{lemma}
\begin{pf}
Fix $p \in\cP$ and $x \in\real$, and suppose that $a(x,Z(x,p,k))
\neq0$. Then there is an open interval containing $x$ on which
$a(\cdot,Z(\cdot,p,k))$ does not vanish and
$b(\cdot,Z(\cdot,p,k))/a(\cdot, Z(\cdot,p,k))$ is continuously
differentiable. Therefore $\ln p$ is, perhaps in a~smaller
neighborhood of $x$, the only solution to the equation~(\ref{eqgenode}) whose derivatives up to order $k$ at $x$ are given
by the components of the vector $Z(x,p,k)$. On the other hand, by
Lemma~\ref{lerichness}(a) there exists $q \in\cP$ such that
$Z(x,q,k) = Z(x,p,k)$ but $z_{k+1}(x,q) \neq z_{k+1}(x,p)$. By
assumption this $q$ is a~solution of~(\ref{eqgenode}), too, with the
same initial conditions. This contradiction to uniqueness is resolved
only if $a(x,Z(x,p,k)) = 0$. Since $p \in\cP$ and $x \in\real$ were
arbitrary, the proof of the lemma is complete.
\end{pf}

Let us combine these facts. Absorbing the (universal) constant $c_p =
c$ in~(\ref{eqexpl1}) into the function $b$, we see that every $p \in
\cP$ satisfies a~differential equation of the form~(\ref{eqgenode})
with $a~= \pd_{22}^2 \mys$. Therefore $\pd_{22}^2 \mys(x,Z(x,p,2)) =
0$ for all $p \in\cP$ and\vadjust{\goodbreak} $x \in\real$ by the reduction principle.
The proof of~(\ref{eqvanish}) is completed on noting that for any $x
\in\real$, $p \in\cP$, $|t| \leq|z_1(x,p)|$ there is a~$q \in\cP$
such that $Z(x,q,2) = (z_0(x,p),t,z_2(x,p))$, by Lemma
\ref{lerichness}(a).

\subsection{Linear dependence on $z_0$} \label{secz0}

The fact that a~local proper score $\mys\in\cR_2$ can be represented
by means of a~$z_2$-independent kernel $K \in\cR_1$ will now be
utilized to show that both $\mys$ and $K$ depend linearly on the
logarithmic score, $z_0$.
%
\begin{proposition}[(Linearity in $z_0$)] \label{propz0lin}
Let $K \in\cR_1$ be the kernel constructed in Section
\ref{secconstruction} from a~given local proper score $\mys\in
\cR_2$. Then $K$ is of the form $K(x,z_0,z_1) = c z_0 +
K_0(x,z_1)$ where $c$ is a~real constant, and $\mys$ is of the
form~(\ref{eqs}) with the same $c$.
\end{proposition}
\begin{pf}
We already know that the score $\mys$ can be represented as in
(\ref{eqtconstrp}), with $K$ not depending on $z_2$. Furthermore, by
the Euler equation~(\ref{eqeeqn}) and Proposition~\ref{propcconst}
there is
some constant $c$ such that
%
\begin{equation} \label{eqeeqn2}
\pd_0 \mys- c = \frac{1}{p} \,\frac{d}{dx} \biggl(p \pd_1 \mys-
\frac{d}{dx} [p \pd_2 \mys] \biggr).
\end{equation}
Using these facts along with the commutation relation $\pd_1\,
\frac{d}{dx} = \frac{d}{dx} \pd_1 + \pd_0$ [cf.~(\ref{eqcommrel2})],
we get
\begin{eqnarray*}
\pd_1 \mys
& = & \pd_1 K - \pd_1 K - z_1 \pd_{11}^2 K - \pd_1 \biggl( \frac{d}{dx} [ \pd
_{1} K] \biggr) + \pd_{01}^2 K \\
& = & - z_1 \pd_{11}^2 K - \frac{d}{dx} [ \pd_{11}^2 K] - \pd_{01}^2 K
+ \pd_{01}^2 K \\
& = & - z_1 \pd_{11}^2 K - \frac{d}{dx} [ \pd_{11}^2 K].
\end{eqnarray*}
On the other hand, we have
\[
\pd_2 \mys
= - \pd_2 \,\frac{d}{dx} \pd_1 K
= - \pd_2 ( \pd_{x1}^2 K + \pd_{01}^2 K \cdot z_1 + \pd_{11}^2 K \cdot
z_2 )
= - \pd_{11}^2 K,
\]
and hence
\[
p^{-1}\,\frac{d}{dx} [ p \pd_2 \mys]
= - z_1 \pd_{11}^2 K - \frac{d}{dx} [\pd_{11}^2 K].
\]
Thus $p \pd_1 \mys= \frac{d}{dx} [ p \pd_2 \mys]$, the
right-hand side of~(\ref{eqeeqn2}) vanishes, and $\pd_0 \mys$ is
constant, $\pd_0 \mys= c$. It follows that $\mys= c z_0 +
g(x,z_1,z_2)$ for some function $g$ independent of~$z_0$, and it
remains to verify the particular forms of $K$ and~$\mys$.

By~(\ref{eqKdef}) and the special form of $\mys$ we have $K - cz_0 =
g - z_1 V - \frac{d}{dx} V$ where now $V = \int_0^{z_1} \pd_2
g(x,t,z_2) \,dt$. But $\pd_2 V = 0$, by~(\ref{eqd2V}) and
(\ref{eqvanish}), and clearly \mbox{$\pd_2 (K - c z_0)=0$}, since $K \in
\cR_1$. Therefore $K_0 = g - z_1 V - \frac{d}{dx} V$ does not depend\vadjust{\goodbreak}
on $z_2$, and it also does not depend on $z_0$ (since neither $g$ nor
$V$ depend on $z_0$), so $K_0 = K_0(x,z_1)$. This completes the proof
of the first claim. The tangent construction based on $K = c z_0 +
K_0(x,z_1)$ then implies, upon observing $\pd_0 K - \int(\pd_0 K) p =
0$, that
\begin{eqnarray*}
\mys & = & K - z_1 \pd_1 K - \frac{d}{dx} [ \pd_1 K] \\
& = & c z_0 + K_0 - z_1 \pd_1 K_0 - \pd_{x1}^2 K_0 - z_2 \pd_{11}^2 K_0,
\end{eqnarray*}
which is the desired representation.
\end{pf}

\subsection{\texorpdfstring{Completion of the proof of Theorem \protect\ref{thcharacterization}}{Completion of the proof of Theorem 3.2}}

The tangent construction based on a~concave functional $\Phi_K$ with $K
\in\cR_1$ yields a~proper score, which is of the form~(\ref{eqs}) if
the kernel $K$ is of the form~(\ref{eqK}). This proves part (a).
Part (b) follows from Propositions~\ref{propR1kernel} and
\ref{propz0lin}. Finally, part (c) is immediate from Theorem
\ref{thGR2007}.

\section{Remaining proofs, supplements and examples} \label{secancillary}

\subsection{\texorpdfstring{Proof of Proposition \protect\ref{propconcave}}{Proof of Proposition 3.3}} \label{secconcave}

Initially, suppose that $K \in\cR_1$ does not depend on $y_0$, so that
$K =
K(x,y_1)$, and is concave in $y_1$ for every fixed $x$. Given $p_0,
p_1 \in\cP$ and $t \in[0,1]$, let $p_t = tp_1 + (1-t)p_0$ and put
$\alpha= tp_1/p_t$, pointwise for every $x \in\real$. Then
$p_t'/p_t = \alpha p_1'/p_1 + (1-\alpha) p_0'/p_0$, whence
\[
K ( \cdot, p_t'/p_t ) \geq \alpha
K ( \cdot, p_1'/p_1 ) + (1-\alpha)
K ( \cdot, p_0'/p_0 )
\]
and so
\begin{eqnarray*}
\Phi_K(p_t)
& = & \int K \biggl( x,\frac{p_t'}{p_t}(x) \biggr) p_t(x) \,dx \\
& \geq& \int\biggl[ \alpha(x) K \biggl(x,\frac{p_1'}{p_1}(x) \biggr)
+ \bigl(1-\alpha(x)\bigr) K \biggl( x,\frac{p_0'}{p_0}(x) \biggr) \biggr] p_t(x) \,dx \\
& = & \int K \biggl( x,\frac{p_1'}{p_1}(x) \biggr) tp_1(x) \,dx
+ \int K \biggl( x,\frac{p_0'}{p_0}(x) \biggr) (1-t) p_0(x) \,dx \\
& = & t \Phi_K(p_1) + (1-t) \Phi_K(p_0).
\end{eqnarray*}
The general case follows by the strict concavity of the entropy
functional $p \mapsto- \int p \ln p$. Concerning the claim about
strict propriety, the pathology described in Remark~\ref{re4} does
not occur within the class $\cP$, because all densities $p \in\cP$ are
strictly positive. Thus, the primitive of $p'/p$ exists throughout~$\real$ and equals $\ln p$ up to a~constant, so that $p'/p = q'/q$
implies $p = q$.~%
\end{pf}

\subsection{Local proper scoring rules of order 1} \label{secone}

The representation~(\ref{eqs}) suggests that local proper scores of
exact order $k = 1$ do not exist. In fact, \citet{ParDawLau} show
that there are no key local score functions of odd order. Within our
framework, we can prove the following.\vadjust{\goodbreak}
%
\begin{proposition} \label{proporder1}
Any local score $\mys\in\cR_1$ that is proper relative to $\cP$
is of the form $\mys= c z_0 + k(x)$ for some $c \leq0$.
\end{proposition}
\begin{pf}
Suppose that $\mys\in\cR_1$ is proper. The Euler equation reduces to
\[
\pd_{0}\mys- \frac{1}{p} \,\frac{d}{dx} [p \pd_1
\mys] = \pd_0 \mys+ z_1 \pd_1 \mys- \pd_{x1}^2 \mys-
z_1 \pd_{01}^2 \mys- z_2 \pd_{11}^2 \mys = c_p\qquad (s = s_p)
\]
in this case. Arguing as in Section~\ref{secreduction}, we find that
$c_p = c$ is independent of $p$ and that $\pd_{11}^2 \mys$ vanishes on
$\real^3$. Therefore there are functions $g, h$ depending only on
$x, z_0$ such that $\mys= z_1 g + h$. Plugging this representation
into the Euler equation gives
\[
c
= z_1 \pd_0 g + \pd_0 h + z_1 g - \pd_x g - z_1 \pd_0 g = z_0' g - \pd
_x g + \pd_0 h,
\]
whence $g = 0$ by another application of the reduction principle.
Thus $\pd_0 h = c$, which means that $\mys=c z_0 + k(x)$. Since $-
z_0$ represents the logarithmic score, $\mys$~can be proper only if $c
\leq0$.
\end{pf}

\subsection{Examples} \label{secexamples}

In the subsequent examples, we keep the notation to a~minimum and
suppress arguments whenever possible.
%
\begin{example} \label{expowerseries}
For $n \geq2$ even and $c \leq0$, let $K = c z_0 - z_1^n$. Then $K
\in\cR_1$, the functional $\Phi_K$ is stricly concave on $\cP$, and
the tangent construction of Proposition~\ref{proptconstr} yields the
score
\begin{eqnarray*}
\mys & = & K - z_1 \pd_1 K - \frac{d}{dx} \pd_1 K + \pd_0 K - \int(\pd
_0 K) q \\
& = & c z_0 - z_1^n + n z_1^n + n(n-1)z_1^{n-2} z_2 + c - c \\
& = & c z_0 + (n-1) (z_1^n + n z_1^{n-2} z_2),
\end{eqnarray*}
which is local of order 2 and strictly proper relative to $\cP$.

Conversely, if $\mys$ is as above, let us carry out the construction
of the associated kernel $K$ described in Section
\ref{secconstruction}. We set
\[
V = \int_0^{z_1} \pd_2 \mys(x,z_0,t,z_2) \,dt = n z_1^{n-1}
\]
and then define $K$ as
\begin{eqnarray*}
K & = & \mys- \biggl[ z_1 + \frac{d}{dx} \biggr] V \\
& = & \mys- n z_1^n - n(n-1)z_1^{n-2} z_2 \\
& = & cz_0 + (n-1) (z_1^n + n z_1^{n-2} z_2) - n z_1^n - n(n-1)
z_1^{n-2} z_2 \\
& = & cz_0 - z_1^n.
\end{eqnarray*}
The construction indeed recovers the kernel $K$ from the score $\mys$.
\end{example}
%
%
\begin{example} \label{exrobust}
The special case $K = - z_1^2$ in the previous example gives the
Hyv\"arinen score, $\mys= z_1^2 + 2 z_2$. Being quadratic in the
log-likelihood derivative, $z_1 = p'/p$, and linear in the second
derivative, $z_2 = p''/p - (p'/p)^2$, this score generally is
sensitive to outliers. For example, within the Gaussian shift-scale
family with mean $\mu$ and variance $\sigma^2$, the Hyv\"arinen score
reduces to $\mys= (x-\mu)^2/\sigma^4 - 2/\sigma^2$.

As an alternative, let us consider the kernel $K = - \ln\cosh z_1$,
which grows only linearly as $z_1$ becomes large. The corresponding
score
%
\begin{equation} \label{eqLCS}
\mys= - \ln\cosh z_1 + z_1 \tanh z_1 + z_2 (1- \tanh^2 z_1)
\end{equation}
appears to be more robust, because as $|y| \rightarrow\infty$,
\[
y \tanh y - \ln\cosh y \rightarrow\ln2,
\]
and the factor of $z_2$ tends to zero exponentially, in that $1-
\tanh^2 y \sim4 \exp{(-2|y|)}$. Of course, the log cosh score
(\ref{eqLCS}) is strictly proper relative to $\cP$, since $K$ is
strictly concave.
\end{example}

\section{Data example: Probabilistic weather forecasting} \label{secdata}

The data example in this section illustrates the use of local and
nonlocal scoring rules in an applied forecasting problem.

Weather forecasting has traditionally been viewed as a~deterministic
enterprise that draws on highly sophisticated, numerical models of the
atmosphere. The advent of ensemble prediction systems in the early
1990s marks a~change of paradigms toward probabilistic forecasting
[\citet{Pal02}, \citet{GneRaf05}].
An ensemble prediction system consists of multiple runs of
numerical weather prediction models, which differ in the initial
conditions and/or the mathematical representation of the atmosphere.
As ensemble forecasts are subject to dispersion errors and biases,
some form of statistical postprocessing is required, for a~happy
marriage of mechanistic and statistical modeling.

\citet{WilHam07} and \citet{BroSmi08} review
statistical postprocessing techniques for ensemble weather forecasts.
State-of-the-art methods include the Bayesian model averaging (BMA)
approach developed by \citet{Rafetal05} and \citet{Sloetal07},
\citet{SloGneRaf10}, and the heterogeneous regression, or ensemble model
output statistics (EMOS), technique of \citet{Gneetal05} and
\citet{ThoGne10}. The BMA approach employs a~mixture distribution, where each mixture component is a~parametric
probability density associated with an individual ensemble member,
with the mixture weight reflecting the member's relative contributions
to predictive skill over a~training period. In contrast, the EMOS
predictive distribution is a~single parametric distribution.

For concreteness, consider an ensemble of point forecasts, $f_1,
\ldots, f_k$, for surface temperature, $x$, at a~given time and
location. The goal is to fit predictive distributions that are as
sharp as possible, subject to them being calibrated
[\citet{GneBalRaf07}]. Let $\phi(x; \mu, \sigma^2)$ denote the
normal density with mean $\mu\in\real$ and variance $\sigma^2 > 0$
evaluated at $x \in\real$. The BMA approach of \citet{Rafetal05}
employs Gaussian components with a~linearly bias-corrected mean. The
BMA predictive density for temperature then becomes
\[
q(x | f_1, \ldots, f_k) =
\sum_{i=1}^k w_i \phi(x; a_i + b_i f_i, \sigma^2)
\]
with BMA weights, $w_1, \ldots, w_k$, that are nonnegative and sum to
1, bias parameters $a_1, \ldots, a_k$ and $b_1, \ldots, b_k$, and a~common variance parameter, $\sigma^2$. The EMOS approach of
\citet{Gneetal05} employs a~single Gaussian predictive density, in that
\[
q(x | f_1, \ldots, f_k) =
\phi(x; a~+ b_1 f_1 + \cdots+ b_k f_k, c + d s^2)
\]
with regression parameters $a$ and $b_1, \ldots, b_k$, and spread
parameters $c$ and $d$, where $s^2$ is the variance of the ensemble
values. The EMOS technique thus is more parsimonious, while the BMA
method is more flexible.

Following the original development in \citet{Rafetal05} and
\citet{Gneetal05}, we apply the BMA and EMOS methods to the
five-member University of Washington Mesoscale Ensemble over the North
American Pacific Northwest [\citet{GriMas02}], at a~prediction
horizon of 48 hours. Here we compare the predictive performance of
the BMA and EMOS density forecasts for surface temperature verifying
in the period of 24 April to 30 June 2000, which is the largest period
common to those used by \citet{Rafetal05} and
\citet{Gneetal05}. The predictive models were fitted on trailing training
periods of length 25 days for BMA and length 40 days for EMOS, as
recommended and described in the aforementioned papers. Overall,
there were 23,691 individual forecast cases at individual
meteorological stations and valid times, when aggregated temporally
and spatially over the test period and the Pacific Northwest,
comprising the states of Washington, Oregon and Idaho, and the
southern part of the Canadian province of British Columbia. All
scores reported are averaged over the 23,691 forecast cases.

%
\begin{table}
\caption{Mean logarithmic score (LS), Hyv\"arinen score (HS), log cosh
score (LCS), quadratic score (QS) and spherical score (SphS) for
statistically postprocessed ensemble forecasts of surface
temperature over the North American Pacific Northwest in April--June
2000, using Bayesian model averaging (BMA) and ensemble model output
statistics (EMOS), respectively. See the text for details}
\label{tabresults}
\begin{tabular*}{\tablewidth}{@{\extracolsep{\fill}}lccccc@{}}
\hline
\textbf{Scoring rule}  & \textbf{LS} & \textbf{HS} & \textbf{LCS}
& \textbf{QS} & \textbf{SphS} \\
\hline
BMA & 2.502 & $-0.113$ & $-0.0572$ & $-0.101$ &
$-0.319$ \\
EMOS & 2.486 & $-0.118$ & $-0.0595$ & $-0.103$ & $-0.321$ \\
\hline
\end{tabular*}
\end{table}

In Table~\ref{tabresults} we assess these forecasts, by computing the
mean score under various local proper scoring rules, namely the
logarithmic score (LS), the Hyv\"arinen score (HS) and the log cosh
score (LCS) introduced in~(\ref{eqLCS}). In addition, we
consider two popular nonlocal scores, namely the quadratic score (QS)
and the spherical score (SphS), defined as
\[
\operatorname{QS}(x,q) = \| q \|_2^2 - 2 q(x)
\quad\mbox{and}\quad
\operatorname{SphS}(x,q) = - \frac{q(x)}{\| q \|_2},
\]
respectively, where \mbox{$\| \cdot\|_2$} denotes the $\mathrm{L}_2$-norm.
These scores are strictly proper relative to the class of the
probability measures with square-integrable Lebesgue densities
[\citet{MatWin76}, \citet{GneRaf07}].

Under all scoring rules, the EMOS technique shows a~slightly lower
(i.e., better) mean score than the BMA method. However, the
differences pale when compared to those between the unprocessed
ensemble forecast and the statistically postprocessed density
forecasts. The unprocessed five-member ensemble gives a~discrete
predictive distribution, namely the empirical measure in $f_1, \ldots,
f_5$, to which the above scores do not apply directly. However, we
can compute the mean score for a~smoothed ensemble forecast, which we
take to be normal, with the first two moments identical to those of
the empirical measure. Under this natural approach, the mean scores
for the smoothed ensemble forecast are very high, reaching 21.4 for
the logarithmic score, $1.14 \times10^4$ for the Hyv\"arinen score,
0.230 for the log cosh score, 0.194 for the quadratic score, and
$-0.217$ for the spherical score, thereby attesting to the benefits of
statistical postprocessing.

\section{Discussion} \label{secdiscussion}

A~scoring rule on the real line is local of order $k$ if the score
depends on the predictive density only through its value, and the
values of its derivatives of order up to $k$, at the realizing event.
It is proper if the expected score is minimized whenever the
predictive density coincides with the density underlying the realizing
event. Supplementing the fundamental work in the recent paper by
\citet{ParDawLau}, we have elaborated a~suitable framework for a~formal characterization of the local proper scoring rules in the
particular, but most relevant, case of order $k \leq2$.

A~practically useful characterization depends on the judicious choice
of a~class $\cS$ of scoring functions, and a~class $\cD$ of predictive
densities, within which scores and densities may vary freely.
Involved therein is a~delicate trade-off, in that narrow classes $\cD$
allow for weak assumptions on the members of $\cS$, but have little,
if any, practical relevance. Our choice of $\cS$---the class
$\cR_2$ of scoring functions growing at most polynomially at
infinity---and of $\cD$---the class $\cP$ of densities decaying faster than the
reciprocal of any polynomial, with log-likelihood derivatives growing at most
polynomially---appears to be usefully general and achieving a~reasonable balance. The balance could easily be shifted, for example,
in favor of more heavy-tailed densities, by adapting the polynomial
growth order in $\cS$.

Counterexamples show that proper scoring rules of practical interest,
such as the Hyv\"arinen score~(\ref{eqHS}), may no longer be
strictly proper relative to any class $\cD$ that contains a~convex
family of densities with a~single common zero. It is thus natural to
assume that all densities in $\cD$ are strictly positive on their
common support, $\Omega$, which then is an interval. The case of
finite boundary points, for example, when $\Omega= (0,\infty)$,
appears to be tractable similarly to the case $\Omega= \real$
considered here, and resulting in essentially the same
characterization. It suffices to impose suitable boundary conditions
at $x = 0$ on the classes $\cS$ and $\cD$, guaranteeing the existence
of integrals and causing the boundary terms in
the proof of Lemma~\ref{lemgrad2} to vanish.

With the resurgence of interest in probabilistic forecasting
[\citet{Gne08}], scoring rules for density forecasts are in increasing
demand. In this context, locality is an appealing property, which we
have studied in this work. A~different argument posits that a~scoring
rule for probabilistic forecasts ought to be sensitive to distance, in
the sense that it rewards the assignment of greater mass not just to
exactly the event or value that is observed, but also to nearby events
[\citet{Sta69}, \citet{JosNauWin09}]. While either approach
has appeal, locality and sensitivity to distance appear to be mutually
exclusive properties, and it is not clear which one is more compelling
[\citet{Mas08}, \citet{WinJos08}]. However, in our meteorological data example as
well as in other experience, local and nonlocal proper scoring rules
generally yield comparable results.

In addition to their use in the assessment of predictive performance,
proper scoring rules play major roles in the theory and practice of
estimation [\citet{Daw07}, \citet{GneRaf07}]. A~striking aspect
is that local proper scoring rules of order $k \geq2$ allow for
statistical inference without knowledge of normalization constants
[\citet{ParDawLau}]. Indeed, this was the motivation for the initial
development by \citet{Hyv05}. The example of the log cosh score
(\ref{eqLCS}) shows that local scores can be less nonrobust than one
might expect. These facets suggest exciting opportunities and novel
prospects particularly in complex settings. Undoubtedly, the
pioneering work of Hyv{\"a}rinen (\citeyear{Hyv05}, \citeyear{Hyv07}),
\citet{DawLau05} and \citet{ParDawLau} has laid the groundwork for a~wide
range of promising future work, both theoretically and
methodologically, and including discrete and multivariate settings
[\citet{DawLauPar}, \citet{Ehm11}], where the tangent
approach may continue to be useful and provide new insight.

\section*{Acknowledgments}

The authors are grateful to Philip Dawid, Steffen
Lauritzen and Matthew Parry for helpful discussions, and sharing
manuscripts and presentations (cf. Remark~\ref{re0}), to Peter
J. Huber for alerting us to the counterexample in Remark~\ref{re4},
to Chris Fraley for providing \textsc{R} code used in Section
\ref{secdata}, and to the Editor, Peter B\"uhlmann, an Associate
Editor and a~referee for constructive feedback.
Tilmann Gneiting thanks the Institute
for Frontier Areas of Psychology and Mental Health in Freiburg,
Germany for hospitality and travel support.


%

\printaddresses

\end{document}